\documentclass[12pt]{amsart}
\usepackage{amssymb}
\usepackage{mathptm}
\usepackage{fullpage} 
\newcommand{\CC}{{\mathbb C}}
\newcommand{\QQ}{{\mathbb Q}}
\newcommand{\RR}{{\mathbb R}}
\newcommand{\ZZ}{{\mathbb Z}}
\newcommand{\HH}{{\mathbb H}}
\newcommand{\bdry}{\partial}
\newcommand{\iso}{\cong}

\newcommand{\maps}{\colon\thinspace}

\DeclareMathOperator{\tr}{tr}
\DeclareMathOperator{\fix}{fix}
\DeclareMathOperator{\rank}{rank}
\DeclareMathOperator{\Affine}{Affine}
\newcommand{\CP}{{ {\mathbb C}\, 
    \mbox{\fontfamily{cmr}\fontshape{n}\selectfont P} }}
\newcommand{\PSL}[2]{\mathrm{PSL}_{#1} #2}
\newcommand{\SL}[2]{\mathrm{SL}_{#1} #2}
\newcommand{\spandef}[2]{{  \left\langle  {#1}  \ \left| \   {#2} \right. \right\rangle }}
\newcommand{\mtext}[1]{\quad\mbox{#1}\quad}

\newcommand{\SnapPea}{\texttt{SnapPea}}
\newcommand{\Snap}{\texttt{Snap}}

\newtheoremstyle{plain}{}{}{\slshape}{}{\bfseries}{.}{0.5em}{}
\theoremstyle{plain} 

\swapnumbers
\newtheorem{theorem}{Theorem}[section]

\newtheorem{lemma}[theorem]{Lemma}

\newtheorem{proposition}[theorem]{Proposition}
\newtheorem{claim}[theorem]{Claim}
\newtheorem{question}[theorem]{Question}
\newtheorem{criterion}[theorem]{Criterion}

\newtheorem{virtual_fibration}[theorem]{Thurston's Virtual Fibration
  Conjecture}
\newtheorem*{main_theorem}{\ref{main_thm}. Theorem}

\theoremstyle{definition}
\newtheorem{defn}[theorem]{Definition}

\theoremstyle{remark}
\newtheorem{rmk}[theorem]{Remark}
\newtheorem{exa}[theorem]{Example}

\makeatletter
  \let\c@theorem=\c@subsection
\makeatother

\def\RCS$#1: #2 ${\expandafter\def\csname RCS#1\endcsname{#2}}
\RCS$Revision: 2.8 $
\RCS$Date: 2001/08/22 14:51:14 $
\renewcommand{\today}{\number\year /\number\month /\number\day}

\begin{document}

\title[Commensurability of 1-cusped 3-manifolds]{Commensurability
   of 1-cusped hyperbolic 3-manifolds}
\author[Calegari]{Danny Calegari}
\address{Department of Mathematics \\ Harvard University \\ Cambridge MA 02138, USA}
\email{dannyc@math.harvard.edu}

\author[Dunfield]{Nathan M. Dunfield}

\address{Department of Mathematics \\ Harvard University \\ Cambridge MA 02138, USA}
\email{nathand@math.harvard.edu}
\date{Version: \RCSRevision, Compile: \today, Last commit: \RCSDate}
\thanks{Both authors were partially supported by the National Science Foundation.}


\begin{abstract} 
  We give examples of non-fibered hyperbolic knot complements in
  homology spheres that are not commensurable to fibered knot
  complements in homology spheres. In fact, we give many examples of
  knot complements in homology spheres where every commensurable knot
  complement in a homology sphere has non-monic Alexander polynomial.
\end{abstract}

\maketitle

\section{Introduction}

For over 20 years, progress in 3-manifold theory has been stimulated
by:
\begin{virtual_fibration}
  Let $M$ be a finite volume hyperbolic 3-man\-i\-fold. Then $M$ has a
  finite cover which is a surface bundle over $S^1$.
\end{virtual_fibration}
Little progress has been made towards the resolution of this
conjecture since it was proposed \cite{Thurston82}.  In fact, there
are few 3-manifolds which do not fiber over $S^1$ but are known to
have finite covers which do fiber (see the references in
\cite[Prob.~3.51]{KirbyList}).  Moreover, Boileau and Wang
\cite{BoileauWang1996} produced infinitely many examples of closed
hyperbolic 3-manifolds for which no solvable cover fibers over $S^1$.
However, fundamental groups of hyperbolic 3-manifolds have huge
numbers of finite index subgroups with a variety of quotients
\cite{Lubotzky95}, and many 3-manifolds fiber over $S^1$ in many
different ways.  Thus it is likely that more complicated classes of
covers provide numerous examples.

A relatively tractable class of 3-manifolds are knot complements in
$S^3$, or more generally, knot complements in rational homology
spheres. If such a 3-manifold fibers over $S^1$, the fiber is a
minimal genus Seifert surface of the knot. Complements of knots in
rational homology spheres rarely cover each other, but much more
frequently they share a common finite cover---that is, they are {\em
  commensurable}. For this class of manifolds, the natural analogue of
Thurston's question is:
\begin{question}
  Let $M$ be a knot complement in a rational homology sphere.  When is
  $M$ commensurable with a fibered knot complement in a rational
  homology sphere?
\end{question}

Here, we give conditions on a knot complement which ensure that it is
not commensurable with a fibered knot complement.  These conditions
are satisfied in many examples, including the complements of a large
number of 2-bridge knots in $S^3$.  These conditions are somewhat
subtle, but we give examples showing that the subtleties are
essential.

Recall the basic:
\begin{criterion}\label{criterion}
  Let $M$ be a knot complement in a rational homology sphere. If the
  Alexander polynomial $\Delta_M$ of $M$ is not monic, then $M$ does
  not fiber over $S^1$.
\end{criterion}
\noindent
Our main idea is to combine this criterion with the fact that the
roots of $\Delta_M$ are related to eigenvalues of reducible
$\PSL{2}{\CC}$-representations of $\pi_1(M)$.

We'll now state the main result.  A 1-cusped hyperbolic 3-manifold $M$
is \emph{generic} if it is not arithmetic and its commensurator
orbifold has a flexible cusp.  The latter condition holds if the cusp
shape of $M$ is not in $\QQ(i)$ or $\QQ(\sqrt{-3})$. Our condition
concerns reducible representations $\rho \maps \pi_1(M) \to
\PSL{2}{\CC}$ whose image is non-abelian.  A representation $\rho \maps
\pi_1(M) \to \PSL{2}{\CC}$ is called \emph{integral} if the trace of
$\rho(\gamma)$ is an algebraic integer for all $\gamma \in \pi_1(M)$.
Otherwise, $\rho$ is \emph{non-integral}. As discussed in
Section~\ref{reducible_and_alex}, whether $M$ has a non-abelian
reducible representation which is non-integral is closely related to
whether $\Delta_M$ has a non-integral root.  When $M$ fibers, its
Alexander polynomial $\Delta_M$ is \emph{monic}, that is, has lead
coefficient $\pm 1$, and so all the roots are integral.  Using the
connection between $\Delta_M$ and reducible representations, it is
easy to show that every non-abelian reducible representation of
$\pi_1(M)$ is integral.

Let $X_0(M)$ denote the geometric component of the
$\PSL{2}{\CC}$-character variety of $\pi_1(M)$ (see
Section~\ref{character_varieties}).  We will show:
\begin{main_theorem} 
  Let $M_1$ be a generic hyperbolic knot complement in a
  $\ZZ/2\ZZ$-homology sphere.  Suppose that the geometric component
  $X_0(M_1)$ contains the character of a non-integral reducible
  representation.  Then $M_1$ is not commensurable to a fibered knot
  complement in a $\ZZ/2\ZZ$-homology sphere.
\end{main_theorem}

Suppose $M_1$ and $M_2$ are commensurable manifolds.  Given a
representation $\rho_1 \maps \pi_1(M_1) \to \PSL{2}{\CC}$ one cannot in
general induce a representation $\rho_2 \maps \pi_1(M_2) \to
\PSL{2}{\CC}$ which is \emph{compatible}, that is, agrees with $\rho_1$
on the fundamental group of the common cover of $M_1$ and $M_2$.
However, when $M_1$ and $M_2$ are generic commensurable 1-cusped
hyperbolic 3-manifolds one knows more.  Here, the existence of a
commensurator \cite{Borel81}, together with Thurston's Hyperbolic Dehn
Surgery Theorem \cite{Thurston82}, gives a natural birational
isomorphism between the geometric components of the character
varieties of $M_1$ and $M_2$ (this is due to the first author
(unpublished) and \cite{LongReid00}).  Moreover, representations of
$\pi_1(M_1)$ coming from its geometric component $X_0(M_1)$ correspond
to compatible representations of $\pi_1(M_2)$ coming from its
geometric component $X_0(M_2)$.  The key to Theorem~\ref{main_thm} is
showing that for a reducible representation of $\pi_1(M_1)$ in
$X_0(M_1)$, the corresponding compatible representation of
$\pi_1(M_2)$ is also reducible.  Then if $X_0(M_1)$ contains the
character of a non-integral reducible representation, there is a
corresponding reducible representation of $\pi_1(M_2)$.  This
representation has to be non-integral as well, and so $M_2$ can't be
fibered.

We end this section with an outline of the rest of the paper.  In
Section~\ref{commensurablity}, we give basic topological restrictions
on when fibered and non-fibered 1-cusped manifolds can be
commensurable.  We also provide constructions of pairs of
commensurable 1-cusped manifolds satisfying these restrictions.
Section~\ref{character_varieties} contains background material about
character varieties.  Section~\ref{rep_comm} discusses representations
of commensurable 3-manifolds.  Section~\ref{reducible_and_alex}
discusses the Alexander polynomial and its connection to reducible
representations.  Section~\ref{pf_of_main_thm} is devoted to the proof
of Theorem~\ref{main_thm}.

In Section~\ref{ex_where_applies}, we show that Theorem~\ref{main_thm}
applies to the complements of many 2-bridge knots, in particular,
to all non-fibered 2-bridge knots $K(p,q)$ where $q < p < 40$.  

In Section~\ref{fiber_and_non}, we give examples of pairs of
commensurable 1-cusped hyperbolic 3-man\-i\-folds exactly one of which
fibers.  These illustrate the necessity and
subtlety of the hypotheses of Theorem~\ref{main_thm}. 

\subsection*{Acknowledgements}  

The first author was partially supported by an NSF Graduate Fellowship
and a NSF VIGRE grant.  The second author was partially supported by
an NSF Postdoctoral Fellowship.  We also thank the referee for helpful
comments, especially those clarifying the definition of an integral
representation.

\section{Commensurability of 1-cusped 3-manifolds}\label{commensurablity}

In this section, we'll discuss some basic necessary conditions for a 1-cusped
non-fibered 3-manifold to be commensurable with a 1-cusped fibered
3-manifold.  We'll also describe some constructions of
pairs of commensurable manifolds which show that these necessary
conditions are satisfied in many examples.

Let $M$ be a 3-manifold.  Given $A$ in $H^1(M; \ZZ)$ we can think of it
in several ways: as a homomorphism of $\pi_1(M)$ to $\ZZ$, as a
homotopy class of maps from $M$ to $S^1$, or as a surface representing
a class in $H_2(M, \bdry M)$ (via Lefschetz duality).  Thought of as a
class of maps from $M$ to $S^1$, it makes sense to say that $A$ is
\emph{representable by a fibration over $S^1$.} 

We'll begin with the question: Suppose $M$ is a 3-manifold which does
not fiber over $S^1$; when does a finite cover of $M$ fiber?  The
fundamental fact here is:
\begin{lemma}[Stallings]\label{Stallings}
  Suppose $M$ is an orientable 3-manifold which does not fiber over
  $S^1$.  Let $N$ be a finite cover which does fiber over $S^1$, with
  $\phi \maps N \to M$ the covering map. If $A$ is a class in $H^1(M;
  \ZZ)$ then the pullback $\phi^*(A)$ in $H^1(N; \ZZ)$ cannot represent
  a fibration over $S^1$. In particular,
  \[
  \rank H^1(N) > \rank H^1(M).
  \] 
\end{lemma}
\begin{proof}
Let $A$ be in $H^1(M, \ZZ)$.  Stallings showed that $A$ can be
represented by a fibration over $S^1$ if and only if the kernel of the
map $A \maps \pi_1(M) \to \ZZ$ is finitely generated
\cite{Stallings62}. So as $M$ does not fiber, the kernel of the $A
\maps \pi_1(M) \to \ZZ$ is not finitely generated.  As $\pi_1(N)$ has
finite index in $\pi_1(M)$, it follows that the kernel of the
restricted map $\phi^*(A) \maps \pi_1(N) \to \ZZ$ is also not finitely
generated.  So $\phi^*(A)$ cannot represent a fibration.
\end{proof}

Now suppose two manifolds $M_1$ and $M_2$ are \emph{commensurable},
that is, they have a common finite sheeted cover $N$.  The following
theorem gives a restriction on when a non-fibered 1-cusped 3-manifold
can be commensurable to a fibered one:

\begin{theorem}\label{cusps_go_up}
  Let $M_1$ and $M_2$ be two commensurable hyperbolic knot complements
  in rational homology spheres.  Suppose $M_2$ fibers over $S^1$
  but $M_1$ does not.  Then a common regular cover must have
  at least $2$ boundary components.
\end{theorem}

\begin{proof}
First, let's make some basic observations.  Throughout, all
(co)homology will have coefficients in $\ZZ$.  Suppose $\phi \maps N
\to M$ is a \emph{regular} finite cover of 3-manifolds.  Let $G =
\pi_1(M)/\pi_1(N)$ be the covering group.  The homomorphism $\phi^*
\maps H^1(M) \to H^1(N)$ is injective, and $\phi^*(H^1(M))$ is exactly
the $G$-invariant cohomology.
  
Now let's prove the theorem.  Suppose $N$ is a common regular cover of
$M_1$ and $M_2$.  Let $\phi_i \maps N \to M_i$ be the covering maps,
and $G_i = \pi_1(M_i)/\pi_1(N)$ the covering groups.  By Mostow
rigidity, we can assume that the covering groups $G_i$ act via
isometries of some fixed hyperbolic metric on $N$.  As the isometry
group of $N$ is finite, so is the group $G = \langle G_1, G_2
\rangle$.
  
From now on, assume that $N$ has only one cusp.  We will show there is
a non-zero $G$-invariant class $A$ in $H^1(N)$.  This gives a
contradiction for the following reason.  Every non-zero class in
$\phi_2^* (H^1(M_2))$ can be represented by a fibration, while by
Lemma~\ref{Stallings} no non-zero class in $\phi_1^* (H^1(M_1))$ can
be represented by a fibration.  But because $A$ is $G$-invariant, $A$
is in $\phi_i^* (H^1(M_i))$ for both $i$, which is impossible.
  
Let $S$ be a Seifert surface for $M_2$ which is a fiber, and let $F =
\phi_2^{-1}(S)$ be the lift to $N$.  The surface $F$ represents a
non-trivial class in $H_2(N, \bdry N)$.  Moreover, since $[\bdry S]$
in $H_1(\bdry M_2)$ is nontrivial, so is $[\partial F]$ is in
$H_1(\partial N)$.  Look at the the class in $H_2(N, \bdry N)$ which
is
\[
  X = \sum_{g \in G} g_*([F]).
\]
  
Consider the restricted coverings $\phi_i \maps \bdry N \to \bdry
M_i$.  Now the covering group $G_i$ acts freely on the torus $\bdry
N$.  Hence $\phi_i$ induces a rational isomorphism on $H_1$, and $G_i$
acts identically on $H_1(\bdry N; \ZZ)$.  Therefore $G$ acts
identically on $H_1(\bdry N; \ZZ)$.  Thus
\[
  \bdry X = \sum_{g \in G} g_*([\bdry F]) = \left| G \right| \cdot [\bdry F] \neq 0.
\] 
So $X$ is non-zero.  If $A$ is the dual class in $H^1(N)$, then $A$ is
the non-zero $G$-invariant class we sought.
\end{proof}

\begin{rmk}
  Suppose $M_1$ and $M_2$ have a common manifold quotient $M$ with a
  torus cusp.  In this case, there are Seifert surfaces generating
  $H_2(M_i,\partial M_i)$ which are pullbacks of a single Seifert
  surface in $M$. In particular, if $M_2$ fibers but $M_1$ does not, the
  commensurator cannot be a manifold, but must be an orbifold.
\end{rmk}

\begin{exa}
Suppose $M_1$ and $M_2$ have a common finite regular cover $N$ with
covering deck groups $G_1$ and $G_2$.  If the $M_i$ are not
hyperbolic, it is not always true that there are actions $G_i'$
isotopic to $G_i$ so that $G=\langle G_1',G_2' \rangle$ is finite.
For example, let $M_1$ be the unit tangent bundle of the $(2,4,4)$
Euclidean triangle orbifold, and $M_2$ the unit tangent bundle of the
$(2,3,6)$ Euclidean triangle orbifold.  These manifolds have $T^3$ as
a common regular cover, with deck groups $G_1 \cong \ZZ/4\ZZ$ and $G_2
\cong \ZZ/6\ZZ$ respectively.  But the action of $\langle G_1,G_2
\rangle$ on $H_1(T^3;\ZZ)$ generates a group isomorphic to $SL(2,\ZZ)$.
Thus, we can't isotope the $G_i$ so that together they generate a
finite group.
\end{exa}

In light of Theorem~\ref{cusps_go_up}, it is worth producing examples
of commensurable 1-cusped manifolds whose common covers have multiple
cusps.

\begin{exa}
Let $N_0 = T^2 \times I$. There are two orientation-preserving
involutions $\phi_1,\phi_2$ of $N_0$ defined by
\[
\phi_1(x,y,t) = (-x,y,1-t), \; \phi_2(x,y,t) = (x,-y,1-t)
\] 
where $(x,y)$ are angular coordinates on $T^2$ and $t$ is the
coordinate on $I = [0,1]$.

Now, the union of the fixed point sets of $\phi_1$ and $\phi_2$ is a
graph $\Gamma \subset T^2 \times 1/2$. A regular neighborhood
$N(\Gamma)$ of $\Gamma$ is a genus $5$ handlebody. Let $N_1 = N_0 -
N(\Gamma)$, and let $N_2$ be the double of $N_1$. Then $N_2$ has
two sets of cusps, $\{C_1,C_2\}$ and $\{D_1,D_2\}$, where $C_1$ and
$C_2$ are the original pair of cusps from $N_0$. The involutions
$\phi_1$ and $\phi_2$ extend to fixed-point--free involutions of $N_2$
which interchange $C_1$ with $C_2$ and $D_1$ with $D_2$.  Notice that
\[
\phi_2 \phi_1^{-1}: T^2 \times 0 \to T^2 \times 0
\] 
is just the involution $(x,y) \to (-x,-y)$. It follows that these
involutions descend to the manifold obtained by equivariantly Dehn
filling $D_1$ and $D_2$.  Let $N_3$ be obtained from $N_2$ by Dehn
filling on both $D_1$ and $D_2$, so that $\phi_1$ and $\phi_2$ both
act on $N_3$. We can also equivariantly surger $N_3$ along some
collection of links to destroy any ``accidental'' additional symmetry
to get $N_3'$, so that the quotients $M_1 = N_3'/\phi_1$ and 
$M_2 = N_3'/\phi_2$ are not homeomorphic, and have no common regular
cover with fewer than two cusps. This equivariant surgery
can also be used to kill off rational homology and ensure $N_3$ is
hyperbolic, so that $M_1$ and $M_2$ can be chosen to be knot
complements in rational homology spheres.
\end{exa}

These examples have the property that for $S$ a Seifert surface in
$M_2$, the class $[\phi_1\phi_2^{-1}(S)]$ is trivial in
$H_2(M_1,\partial M_1)$.  So the proof of Theorem~\ref{cusps_go_up}
does not apply here, and we cannot conclude anything about whether
$M_1$ and $M_2$ mutually fiber or do not fiber over $S^1$.

These examples cannot be chosen to be knot complements in
$\ZZ/2\ZZ$-homology spheres because of the very existence of a $2$-fold
cover. A modification of this construction gets around this
difficulty.

\begin{exa}
Let $L$ be the unlink in $S^3$ on $n$ components. Arrange these
symmetrically so that there is a rotation $r$ with axis $\alpha$
permuting the components of $L$. Let $s$ be a rotation of order $n$
fixing each component, translating $\alpha$, and fixing another axis
$\beta$ which links each component of $L$. Let $\phi_1 = rs$ and
$\phi_2 = rs^k$ for some $k>1$ coprime with $n$. Then
$M_i=(S^3-L)/\phi_i$ is a knot complement in a lens space which is a
$\ZZ/2\ZZ$-homology sphere for $n$ odd. By equivariant surgery, we can
make the $M_i$ hyperbolic knot complements in $\ZZ/2\ZZ$-homology spheres
whose smallest common cover has $n$ cusps.
\end{exa}

\section{Character varieties}\label{character_varieties}

Here, we review the part of the theory of character varieties of
3-manifolds that we will need for Theorem~\ref{main_thm}.  For
details, see \cite{CullerShalen83, ShalenHandbook}.  For the
technicalities of the $\PSL{2}{\CC}$ as opposed to $\SL{2}{\CC}$ case see
\cite{BoyerZhang98}.

\begin{defn}
For $M$ a compact 3-manifold, we let
\[
R(M) = \text{Hom}(\pi_1(M),\PSL{2}{\CC})
\] 
denote the $\PSL{2}{\CC}$ representation variety of $M$.

Further, let
\[
X(M) = \bigcup_{\rho \in R(M)} \tr_\rho^2 \subseteq \CC^{\pi_1(M)}
\] 
denote the $\PSL{2}{\CC}$ character variety of $M$.
\end{defn}

Since $\pi_1(M)$ is finitely generated, $R(M)$ is an affine algebraic
variety over $\CC$. For each $\gamma \in \pi_1(M)$, the function
$\tau_\gamma:R(M) \to \CC$ defined by
\[
\tau_\gamma(\rho)=\tr_\rho^2(\gamma)
\] 
is a regular function on $R(M)$. As $\pi_1(M)$ is finitely generated,
the ring of functions generated by the $\tau_\gamma$ is finitely
generated. It follows that $X(M)$ in $\CC^{\pi_1(M)}$ projects
isomorphically to an algebraic subvariety of $\CC^G$ for some finite
subset $G \subset \pi_1(M)$.  Therefore $X(M)$ has the structure of an
affine algebraic variety over $\CC$.

Away from the reducible locus, the action of $\PSL{2}{\CC}$ on $R(M)$
by conjugation is transverse, and the natural projection
$R(M)/\PSL{2}{\CC} \to X(M)$ is injective on a Zariski open set.  Let
$t \maps R(M) \to X(M)$ be the projection.  For a character $\chi$ in
$X(M)$, $t^{-1}(\chi)$ either consists solely of the conjugates of a
single irreducible representation, or $t^{-1}(\chi)$ consists of
reducible representations.  In the latter case, the reducible
representations in $t^{-1}(\chi)$ may not all be conjugate (it is easy
to see that the closure of the orbit of a non-abelian reducible
representation contains an abelian reducible representation).

We will need the following, which is Proposition~6.2 of \cite{CGLS}.
\begin{lemma}\label{non_abel_red}
  Let $X$ be an irreducible component of $X(M)$ which contains the
  character of an irreducible representation.  Let $\chi \in X$ be the
  character of a reducible representation.  Then there is a
  non-abelian reducible representation with character $\chi$.
\end{lemma}
The lemma follows from the fact that the fibers of $t$ over $X$ are
all at least 3-dimensional whereas the orbit under conjugation of an
abelian reducible representation is 2-dimensional.

If $M$ is an orientable finite volume hyperbolic 3-manifold, there is
a unique discrete faithful representation
\[
\rho_\delta:\pi_1(M) \to \PSL{2}{\CC}
\] 
up to conjugacy in $O(3,1) \cong \PSL{2}{\CC} \rtimes \ZZ/2\ZZ$.  Up to
conjugacy in $\PSL{2}{\CC}$, there are a pair of such representations
which differ by complex conjugation, and their characters may occur in
different irreducible components of $X(M)$. Fixing an orientation of $M$
fixes a $\PSL{2}{\CC}$-conjugacy class of discrete faithful
representations.  We will assume our manifolds are oriented in what
follows, and therefore that $\rho_\delta$ is well-defined up to
conjugacy in $\PSL{2}{\CC}$.

\begin{defn}
  Let $M$ be a finite volume hyperbolic 3-manifold. Let $X_0(M)$
  denote the irreducible component of $X(M)$ containing the character
  of the discrete faithful representation $\rho_\delta$.  The
  component $X_0(M)$ is called the \emph{geometric component} of
  $X(M)$.
\end{defn}

\section{Representations of commensurable manifolds}\label{rep_comm}

Let $M_1$ and $M_2$ be commensurable hyperbolic 3-manifolds with
common finite cover $N$.  Two representations $\rho_i \maps \pi_1(M_i)
\to \PSL{2}{\CC}$ are said to be \emph{compatible} if they agree on
$\pi_1(N)$.  For instance, if $\rho_1$ is a discrete faithful
representation for $M_1$, then Mostow rigidity implies that there is a
discrete faithful representation $\rho_2$ of $\pi_1(M_2)$ which is
compatible with $\rho_1$.  The property of having compatible
representations extends to other representations whose characters are
in $X_0(M_1)$.

\begin{proposition}\label{compat_on_geom_comp}
  Suppose $M_1$ and $M_2$ are generic commensurable orientable hyperbolic
  3-manifolds with one cusp.  Let $\chi_1$ be a character in the
  geometric component $X_0(M_1)$.  Then there exist compatible
  representations $\rho_i$ of $\pi_1(M_i)$ such that $\rho_1$ has
  character $\chi_1$ and the character of $\chi_2$ lies in a geometric
  component $X_0(M_2)$.
\end{proposition}

The reason for weaseling around with $\chi_1$ instead of just starting
with $\rho_1$ is that for characters of reducible representations,
there can be different conjugacy classes of representations with that
character.  It turns out that you can't always specify $\rho_1$, but
only $\chi_1$, in this case.

\begin{proof}
As the $M_i$ are non-arithmetic, they cover a common orientable
commensurator orbifold $Q$ \cite{Borel81}. Let $p_i \maps M_i \to Q$
be the (orbifold) covering maps.  The inclusion of $\pi_1( M_i)$ into
$\pi_1(Q)$ induces a map ${p_i}_* \maps X_0(Q) \to X_0(M_i)$ via
restriction of representations.  Because the $M_i$ are generic, $Q$
has a flexible cusp, and the variety $X_0(Q)$ is also a complex curve.
In fact, ${p_i}_*$ is a birational isomorphism, though we will not
need this \cite{LongReid00}.  The main step is:

\begin{lemma}\label{surjects} 
The map ${p_1}_* \maps X_0(Q) \to X_0(M_1)$ is onto.
\end{lemma}
  
\begin{proof}[Proof of Lemma.]
The map ${p_1}_*$ is a non-constant map of irreducible affine
algebraic curves over $\CC$.  Let $\bar{X}_0(Q)$ denote the smooth
projective model of $X_0(Q)$.  The curve $\bar{X}_0(Q)$ is the
normalization of $X_0(Q)$ compactified by adding an ideal point for
each end of $X_0(Q)$ \cite[\S 1.5]{CGLS}.  Similarly, let
$\bar{X}_0(M_1)$ be the smooth projective model of $X_0(M_1)$.  The
map ${p_1}_*$ induces a regular map of the same name between
$\bar{X}_0(Q)$ and $\bar{X}_0(M_1)$ (this map is just a branched
covering of closed Riemann surfaces).  Let $\chi_1$ be a point in
$\bar{X}_0(M_1)$ which corresponds to a character---that is, not an
ideal point.  As the map from $\bar{X}_0(Q)$ to $\bar{X}_0(M_1)$ is
surjective, choose $\chi_0$ in $\bar{X}_0(Q)$ with ${p_1}_*(\chi_0) =
\chi_1$.  We need to show that $\chi_0$ is not an ideal point.
Suppose that $\chi_0$ is an ideal point.  By Proposition~1.4.4 of
\cite{CullerShalen83} there is some $\gamma$ in $\pi_1(Q)$ for which
$\tr^2_\gamma(\chi_0) = \infty$.  That is, there is some element of
$\pi_1(Q)$ which acts by a hyperbolic isometry on the simplicial tree
associated to the ideal point $\chi_0$.  Now for any $n>0$, $\gamma^n$
also acts by a hyperbolic isometry on the tree and so
$\tr^2_{\gamma^n}(\chi_0) = \infty$.  As $\pi_1(M_1)$ is of finite
index in $\pi_1(Q)$, we can choose $n$ so that $\gamma^n$ is in
$\pi_1(M_1)$.  But then $\tr^2_{\gamma^n} \chi_0 = \tr^2_{\gamma^n}
\chi_1 = \infty$, contradicting that $\chi_1$ is the character of a
representation.  So $\chi_0$ is not an ideal point and hence ${p_1}_*
\maps X_0(Q) \to X_0(M_1)$ is onto.
\end{proof}
  
Now to finish the proof of the theorem, let $\chi_1 \in X_0(M_1)$.  By
the lemma, there is some character $\chi_0$ in $X_0(Q)$ with
${p_1}_*(\chi_0) = \chi_1$.  Let $\rho_0$ be a representation with
character $\chi_0$.  Then the restrictions of $\rho_0$ to the
subgroups $\pi_1(M_i)$ give a pair of compatible representations with
the required properties.
\end{proof}

\section{The Alexander polynomial and reducible representations}
\label{reducible_and_alex}

Let $M$ be a knot complement in a rational homology sphere.  Let $N$
denote the universal free abelian cover of $M$.  That is, set $H =
H_1(M;\ZZ)/(\mathrm{torsion})$ and take the covering corresponding to
the kernel of the natural homomorphism $\pi_1(M) \to H$. Then $N$ is a
regular covering of $M$, and the group $\pi_1(M)/\pi_1(N) = H \iso \ZZ$
acts on $N$ by deck transformations. It follows that $H_1(N;\ZZ)$ has
the natural structure of a $\ZZ[H]$-module.  If $t$ denotes the
generator of $H$, then $H_1(N;\ZZ)$ is a $\ZZ[t,t^{-1}]$-module.  The
\emph{Alexander polynomial} $\Delta_M$ of $M$ is an invariant of this
module.  In the case that $H_1(N;\ZZ)$ is cyclic, that is $H_1(N;\ZZ) =
\ZZ[t,t^{-1}]/p(t)$, the polynomial $\Delta_M$ is just $p(t)$.  In
general, $\Delta_M$ is the greatest common divisor of the $0$-th
elementary ideal of a matrix which presents $H_1(N;\ZZ)$ as a
$\ZZ[t,t^{-1}]$-module.  The Laurent polynomial $\Delta_M$ is only
defined up to multiplication by a unit $\pm t^n$ in $\ZZ[t,t^{-1}]$.  A
key property for us is that the Alexander polynomial is symmetric,
that is $\Delta_M(t^{-1}) = \pm t^n \Delta_M(t)$ \cite{Turaev75}.  For
more on the Alexander polynomial see \cite{Rolfsen76}.

If $M$ is a surface bundle over $S^1$ with fiber $F$ and
monodromy $\phi:F \to F$, then $N = F \times \RR$ and the action of $t$
on $H_1(N;\ZZ)$ is exactly equal to the action of $\phi_*$ on
$H_1(F;\ZZ) = H_1(N;\ZZ)$.  In this case, $\Delta_M$ is the
characteristic polynomial of $\phi_*$. Since $\phi$ is a
homeomorphism, $\phi_*$ is an automorphism, and $\Delta_M$ is
{\em monic.}

Now consider a non-abelian reducible representation $\rho \maps \pi_1(M) \to
\SL{2}{\CC}$. Conjugate $\rho$ so that its image is upper-triangular.
Given $\gamma \in \pi_1(M)$ we will say that the \emph{primary
  eigenvalue} of $\rho(\gamma)$ is the $(1,1)$ entry of the matrix
$\rho(\gamma)$.  This is well defined for the following reason.  Since
$\rho$ has non-abelian image, there is a unique line $L$ in $\CC^2$
which is left invariant by all $\rho(\pi_1(M))$.  The primary eigenvalue
of $\rho(\gamma)$ is just the eigenvalue of $\rho(\gamma)$ with
eigenspace $L$.

A reducible representation into $\SL{2}{\CC}$ has meta-abelian
(two-step solvable) image.  The Alexander polynomial of $M$ is an
invariant of the maximal meta-abelian quotient of $\pi_1(M)$, so it's
not surprising that it is related to non-abelian reducible
representations of $\pi_1(M)$ into $\SL{2}{\CC}$.  For knots in
$\ZZ$-homology spheres, the statement is:

\begin{theorem}[de Rham]\label{de_Rham}
  Let $M$ be a knot complement in a $\ZZ$-homology sphere. Let $\mu$ in
  $\pi_1(\bdry M)$ be a meridian.  The following are equivalent:
  \begin{itemize}
    \item There is a non-abelian reducible representation
    \[
    \rho \maps \pi_1(M) \to \SL{2}{\CC}
    \] 
    such that $\rho(\mu)$ has primary eigenvalue
    $m$.
  
  \item $m^2$ is a root of $\Delta_M(t)$.
    
\end{itemize}
\end{theorem}

More generally, for knots in $\QQ$-homology spheres there is a similar
connection that is a bit harder to state.  Let $\rho \maps \pi_1(M)
\to \SL{2}{\CC}$ be a non-abelian reducible representation.  Then
$\rho$ acts on $\CP^1$ and has a unique fixed point.  Translating that point
to $\infty$, $\rho$ can be interpreted as a homomorphism from
$\pi_1(M)$ into the (complex) affine group of $\CC$:
\[
   \Affine(\CC) = \left\{ \text{maps } z \mapsto a z + b \right\}
        \iso \left\{ \left(
              \begin{array}{cc} a & b \\ 0&  1 \end{array} \right) \right\}
              \mtext{where $a \in \CC^\times$ and $b \in \CC$.}
\] 
Define a homomorphism $x_\rho \maps \pi_1(M) \to \CC^\times$ by setting
$x_\rho(\gamma) = a$, where $a$ is the homothety part of
$\rho(\gamma)$ thought of as an element of $\Affine(\CC)$. Note that
$x_\rho(\gamma)$ is just the square of the primary eigenvalue of
$\rho(\gamma)$ regarded as being in $\SL{2}{\CC}$.  The map $x_\rho$ is
often called the character of $\rho$ but to prevent confusion we'll
avoid this practice.  The case of knots in $\QQ$-homology spheres is
more complicated than the $\ZZ$ case because not every homomorphism $x
\maps \pi_1(M) \to \CC^\times$ factors through the \emph{free}
abelianization of $\pi_1(M)$.  For those homomorphisms that do,
Theorem~\ref{de_Rham} in this context is just:
\begin{theorem}\label{free_ab}
  Let $M$ be a knot complement in a $\QQ$-homology sphere. Let $H$ be
  the free abelianization of $\pi_1(M)$.  Let $x \maps H \to \CC^\times$ be a
  homomorphism.  Then the following are equivalent:
  \begin{itemize}
  \item There is a non-abelian reducible representation
      \[
      \rho \maps \pi_1(M) \to \SL{2}{\CC}
      \] 
      with $x_\rho = x$.
    \item $x'(\Delta_M) = 0$, where $x'$ is the map $x$ induces from
      $\ZZ[H]$ to $\CC$.
\end{itemize}
\end{theorem}
For a proof, see e.g.~\cite[\S 3]{McMullenNorm}. 

\subsection{Non-integral reducible representations}

From the introduction, recall that a representation $\rho$ of
$\pi_1(M)$ into $\PSL{2}{\CC}$ is \emph{integral} if the trace of every
$\rho(\gamma)$ is an algebraic integer.  Otherwise $\rho$ is
\emph{non-integral}.  Let's reformulate this a little.  Consider a
matrix $A \in \PSL{2}{\CC}$.  An eigenvalue $\lambda$ of $A$
(well-defined up to sign) is a root of the monic polynomial $x^2 \pm
\tr(A) x + 1$.  Therefore if $\tr(A)$ is integral, the eigenvalues of
$A$ are also algebraic integers.  The converse is clearly true, so $A$
has integral trace if and only if \emph{both} eigenvalues of $A$ are
integral.  Now if $\lambda$ is an eigenvalue of $A$, then the other
eigenvalue is $1/\lambda$.  Thus we see that $A$ has integral trace if
and only if the eigenvalues of $A$ are algebraic units.  In the
terminology of the last section, a reducible representation $\rho$ is
integral if and only if the primary eigenvalues of all the
$\rho(\gamma)$ are algebraic units.  Since the primary eigenvalue of
$\rho(\gamma^{-1})$ is the inverse of the primary eigenvalue of
$\rho(\gamma)$, we see that a reducible representation is integral if
and only if every primary eigenvalue is an algebraic integer.  If we
think of $\rho$ as a homomorphism to $\Affine(\CC)$, we see that $\rho$
is integral if and only if the homothety of each $\rho(\gamma)$ is
integral.

For knots in $\ZZ$-homology spheres, the Alexander polynomial
determines the existence of non-abelian reducible representations
which are non-integral:

\begin{proposition}\label{non_int_non_monic}
Let $M$ be a knot complement in a $\ZZ$-homology sphere.  Then $M$ has
a non-abelian reducible representation into $\PSL{2}{\CC}$ which is
non-integral if and only if $\Delta_M$ is not monic.
\end{proposition}

\begin{proof} 
As $M$ is a knot complement in a $\ZZ/2\ZZ$-homology sphere, the
cohomology group $H^2(M; \ZZ/2\ZZ)$ vanishes and every representation
into $\PSL{2}{\CC}$ lifts to $\SL{2}{\CC}$, so we're free to think about
$\SL{2}{\CC}$-representations instead.  Consider a reducible
representation $\rho$.  Since $H_1(M) = \ZZ$ is generated by $\mu$, we
see that the primary eigenvalue of any $\rho(\gamma)$ is a power of
the primary eigenvalue of $\rho(\mu)$.  Thus, it follows from
Theorem~\ref{de_Rham} that $M$ has a non-integral reducible
representation if and only if $\Delta_M$ has a root which is not an
algebraic unit.
  
Suppose that $\Delta_M$ is not monic.  Then $\Delta_M$ has a
non-integral root provided that $\Delta_M$ is not an integer multiple
of a monic integer polynomial.  As we're in the $\ZZ$-homology sphere
case, $\Delta_M(1) = \pm 1$, and this can't happen. So $\Delta_M$ has
a non-integral root, and thus a non-abelian reducible representation
which is non-integral.
  
Now suppose that $\Delta_M$ is monic.  Then all the roots of
$\Delta_M$ are algebraic integers.  Let $\alpha$ be a root of
$\Delta_M$.  Because $\Delta_M$ is symmetric, $1/\alpha$ is also a
root of $\Delta_M$ and so is integral.  Thus all the roots of
$\Delta_M$ are algebraic units.  So all the non-abelian reducible
representations of $\pi_1(M)$ are integral, completing the proof of
the theorem.
\end{proof}

In the general $\ZZ/2\ZZ$-homology sphere case, there isn't an easy
statement like this because Theorem~\ref{free_ab} only applies to
representations coming from certain homomorphisms to $\CC^\times$.  It
is true that if $M$ has a non-integral representation then $\Delta_M$
is non-monic (to prove this, a nice point of view is the theory of BNS
invariants \cite{Dunfield:norms, BieriNeumannStrebel87, Brown87}).
However, if $\Delta_M$ is non-monic, $M$ need not have a non-integral
reducible representation (e.g.~the SnapPea census manifold $m261$).
Nor is the above proposition true for the $\ZZ/2\ZZ$-homology sphere
case if one replaces the non-monic hypothesis with ``has a
non-integral root'' (to see that the ``only if'' direction is false, take
the complement of a fibered knot in $S^3$ connected sum with a lens
space).

Regardless, the following proposition, which in the $\ZZ$-homology
sphere case follows immediately from
Proposition~\ref{non_int_non_monic} and Criterion~\ref{criterion}, is
easy to prove in general.
\begin{lemma}\label{reps_of_fibered}
  Let $M$ be a knot complement in a $\ZZ/2\ZZ$-homology sphere.  If
  $M$ fibers over $S^1$ then every non-abelian reducible
  representation of $\pi_1(M)$ into $\PSL{2}{\CC}$ is integral.
\end{lemma}
\begin{proof}
Let $\rho \maps \pi_1(M) \to \Affine(\CC)$ be a lift of a given
non-abelian reducible $\PSL{2}{\CC}$ representation.  As $M$ fibers
over $S^1$, the universal abelian cover of $M$ is of the form $F
\times \RR$ where $F$ is a compact surface (here $F$ is some finite
abelian cover of a fiber in the fibration of $M$).  As $\pi_1(F)$ is
the commutator subgroup of $\pi_1(M)$, the representation $\rho$ takes
$\pi_1(F)$ to a finitely generated abelian subgroup $G$ consisting of
translations. The subgroup $G$ is non-trivial as $\rho$ is
non-abelian.  For each $\gamma \in \pi_1(M)$ we need to show that if
$A = \rho(\gamma) = (z \mapsto a z + b)$ then the homothety $a$ is an
algebraic integer.  The action of $A$ by conjugation on the normal
subgroup $G$ takes an element $(z \mapsto z + \tau)$ to $(z \mapsto z
+ a \tau )$.  So $A$ is a group automorphism of the lattice $G \iso
\ZZ^n$.  Thought of as an element of $\SL{n}{\ZZ}$, the map $A$
satisfies its characteristic polynomial $f(t)$, which is a monic
polynomial with integer coefficients.  Let $B = (z \mapsto z + \tau)$
be a non-identity element of $G$.  If we act on $B$ via $f(A)$ we get
that $f(a) \tau = 0$.  Thus $f(a) = 0$ and $a$ is an algebraic
integer.  So $\rho$ is integral.
\end{proof}

\section{Inducing reducible representations of commensurable manifolds}
\label{pf_of_main_thm}

This section is devoted to proving:

\begin{theorem}\label{main_thm}
  Let $M_1$ be a generic hyperbolic knot complement in a
  $\ZZ/2\ZZ$-homology sphere.  Suppose that the geometric component
  $X_0(M_1)$ contains the character of a non-integral reducible
  representation.  Then $M_1$ is not commensurable to a fibered knot
  complement in a $\ZZ/2\ZZ$-homology sphere.
\end{theorem}

\begin{proof}  
Suppose that $M_1$ is commensurable to another knot complement $M_2$
in a $\ZZ/2\ZZ$-homol\-ogy sphere.  Call the common finite cover $N$.
We will show that $M_2$ has a non-abelian reducible
representation which is non-integral, and so cannot fiber.
  
Let $\chi_1$ in $X_0(M_1)$ be the character of a non-integral
reducible representation.  As $M_1$ is generic, by
Proposition~\ref{compat_on_geom_comp}, there are representations
$\rho_i \maps \pi_1(M_i) \to \PSL{2}{\CC}$ which agree on $\pi_1(N)$
where the character of $\rho_1$ is equal to $\chi_1$.  In particular,
$\rho_1$ is reducible and non-integral ($\rho_1$ may be abelian,
because we don't get to pick $\rho_1$, just $\chi_1$).  Also, the
character of $\rho_2$ is in $X_0(M_2)$.
  
We will show

\begin{claim}\label{claim}
The representation $\rho_2$ of $\pi_1(M_2)$ is reducible.
\end{claim}
  
Assuming the claim, let's prove that $M_2$ is not fibered.  Pick
$\gamma$ in $\pi_1(M_1)$ such that $\rho_1(\gamma)$ has non-integral
trace.  Then for any $n>0$, the matrix $\rho_1(\gamma^n)$ also has
non-integral trace as its eigenvalues are powers of those of
$\rho_1(\gamma)$.  Since $\pi_1(N)$ is of finite index in
$\pi_1(M_1)$, choose an $n$ such that $\gamma^n$ is in $\pi_1(N)$.
But then $\gamma^n$ is in $\pi_1(M_2)$ as well, and so
$\rho_2(\gamma^n) = \rho_1(\gamma^n)$ has non-integral trace.  Thus
$\rho_2$ is non-integral.  Since the character of $\rho_2$ is in
$X_0(M_2)$, by Lemma~\ref{non_abel_red} there is a non-abelian
reducible representation $\rho'_2$ which has the same character as
$\rho_2$.  As $\rho'_2$ has the same character as $\rho_2$, it is
non-integral.  By Lemma~\ref{reps_of_fibered}, $M_2$ does not fiber
over $S^1$.  This completes the proof of the theorem modulo the claim.
   
Now let's go back and prove Claim~\ref{claim}.  Let $\Gamma =
\pi_1(M_2)$ and $\Gamma' = \pi_1(N)$.  Now $\rho_2$ restricted to
$\Gamma'$ is the same as $\rho_1$ restricted to $\Gamma'$, and
$\rho_1$ is reducible.  Thus $\rho_2$ is reducible on $\Gamma'$.  The
subgroup $\Gamma'$ is of finite index in $\Gamma$, so we can replace
it by a finite index normal subgroup of $\Gamma$.  Let $G =
\rho_2(\Gamma)$ and $G' = \rho_2(\Gamma')$, two subgroups of
$\PSL{2}{\CC}$.  Note that $\Gamma'$ is not the trivial subgroup
because $\rho_1$ is non-trivial, in fact non-integral, on any finite
index subgroup of $\pi_1(M_1)$.
 
Now suppose that $\rho_2$ is irreducible, that is, the fixed point set
of $G$ acting on $\bdry \HH^3$ is empty.  As $G'$ is reducible,
$\fix(G')$ is either 1 or 2 points.  As $G'$ is normal in $G$, the set
$\fix(G')$ is $G$-invariant.  So if $\fix(G')$ consisted of a single
point, $G$ would be reducible as well.  So $\fix(G')$ is 2 points.
Look at the homomorphism $h \maps G \to \ZZ/2\ZZ$ where $\ZZ/2\ZZ$ is
thought of as the symmetric group on $\fix(G')$.  The homomorphism $h$
is non-trivial as $G$ is irreducible.  Any $A \in G$ leaves invariant
the geodesic $g$ joining the two points $\fix(G')$.  If $h(A) = 1$
then $A$ acts on $g$ by an orientation reversing isometry, and $A$ has
order $2$.  Note that $G$ is meta-abelian, as the kernel of $h$ is
abelian because it consists of isometries which fix the pair of points
$\fix(G')$.

To finish the proof of the claim, we look at $H = \rho_2(\pi_1(\bdry
M_2))$.  We claim that $H$ is finite.  Let $\mu_2$ in $\pi_1(\bdry
M_2)$ be a meridian, that is, Dehn filling in along $\mu_2$ yields a
$\ZZ/2\ZZ$-homology sphere.  Let $\lambda_2$ in $\pi_1(\bdry M_2)$ be a
longitude, that is, a generator of the kernel $H_1(\bdry M_2, \ZZ) \to
H_1(M_2, \ZZ)$.  If $M$ were the complement of a knot in a
$\ZZ$-homology sphere, $(\mu_2, \lambda_2)$ would be a basis of
$\pi_1(\bdry M_2)$.  In general, $(\mu_2, \lambda_2)$ generate a
finite index subgroup of $\pi_1(\bdry M_2)$.  As $\mu_2$ generates
$H_1(M_2, \ZZ_2)$, we must have $h \circ \rho_2(\mu_2) = 1$ and
$\rho_2(\mu_2)$ has order two.

We claim that since $M_2$ is a knot complement in a $\ZZ/2\ZZ$-homology
sphere, if $K$ is the kernel of the unique surjection $\pi_1(M_2) \to
\ZZ/2\ZZ$ then $\lambda_2 \in [ K, K ]$.  Consider a Seifert surface $S$
for $M$.  The surface $S$ has $[\bdry S]$ equal to $[\lambda_2]$ in
$H_1(\bdry M)$.  We can explicitly construct the cover $M_2'$
corresponding to $K$ by gluing together two copies of $M_2$ cut along
$S$.  Thus we see that $S$ lifts to $M_2'$. This shows that the
boundary of $S$, namely $[\lambda_2]$, is $0$ in $H_1(M_2', \ZZ)$.  Thus
$\lambda_2$ is in $[K, K]$.  Therefore, as $\rho_2(K)$ is an abelian
group of isometries fixing $\fix(G')$, we have $\rho_2(\lambda_2) =
I$.  So the subgroup of $H$ generated by the images of $(\mu_2,
\lambda_2)$ is finite, in fact has order 2.  Thus $H$ itself is
finite.
  
Now we'll argue that $H$ is infinite, yielding a contradiction.  Look
at $M_1$ and in particular at $\rho_1(\pi_1(\bdry M_1))$.  Let $\mu_1$
be a meridian in $\pi_1(\bdry M_1)$.  As $\rho_1$ is non-integral, it
is easy to see from the homomorphism $x_{\rho_1} \maps \pi_1(M_1) \to
\CC^\times$ that the $\gamma$ in $\pi_1(M_1)$ with non-integral trace
are exactly those $\gamma$ which are non-zero in $H_1(M_1,
\ZZ)/(\mathrm{torsion})$.  Therefore, $\rho(\mu_1)$ has non-integral
trace.  In particular, $\rho_1(\mu_1)$ has infinite order, and hence
$\rho_1(\pi_1(\bdry M_1))$ is infinite.  As $\pi_1(\bdry M_1)$ shares
a finite index subgroup with $\pi_1(\bdry M_2)$, $H$ shares a finite
index subgroup with $\rho_1(\pi_1(\bdry M_1))$.  Thus $H$ is infinite.
But we've already shown that $H$ is finite.  This is a contradiction,
and so $\rho_2$ must be reducible.  This proves Claim~\ref{claim} and
thus the theorem.
\end{proof}

\section{2-bridge knots to which the theorem applies}
\label{ex_where_applies}

Theorem~\ref{main_thm} applies to many 2-bridge knots in $S^3$.  A
2-bridge knot is determined by a pair of relatively prime odd
integers $(p,q)$ with $0 < q < p$ (for background see \cite[\S
12]{BurdeZieschang}, \cite{HatcherThurston}).  In this section, we
describe computations which show:
\begin{theorem}\label{2-bridge}
Let $K(p,q)$ be a 2-bridge knot such that $p < 40$.  Let $M$ be the
exterior of $K(p,q)$.  If $M$ does not fiber over the circle, then $M$
satisfies the hypotheses of Theorem~\ref{main_thm}, and so $M$ is not
commensurable to a fibered knot complement in a $\ZZ/2\ZZ$-homology
sphere.
\end{theorem}

Let $K(p,q)$ be a 2-bridge knot, and $M$ be its exterior.   We will follow
\cite{HildenLozanoMontesinos96}, where Hilden, Lozano, and Montesinos,
building on work of Burde and Riley, give a simple method for computing
the $\PSL{2}{\CC}$-character variety $X(M)$.

The standard presentation of $\pi_1(M)$ has as generators two elements
$a$ and $b$, each of which is a meridian at the top of one of the two
bridges.  As $a$ and $b$ are conjugate, we can take coordinates on
$X(M)$ to be $x = \tr_{a^2}$ and $z = \tr_{ab}$ (the latter makes
sense even in $\PSL{2}{\CC}$ because $a$ and $b$ are conjugate).  Thus
$X(M)$ is a plane curve.  There is a polynomial with integer
coefficients $f(x,z)$ such that $X(M)$ is the set of points in $\CC^2$
satisfying $f(x,z)=0$.  Section~5 of \cite{HildenLozanoMontesinos96}
gives a simple recursive procedure for computing this polynomial.

Let $M$ be the complement of a 2-bridge knot which does not fiber.
Because 2-bridge knots are alternating, this is equivalent to
$\Delta_M$ being non-monic (see e.g.~\cite[\S 13.C]{BurdeZieschang}).
To decide if $M$ satisfies the hypotheses of Theorem~\ref{main_thm},
we first need to factor the polynomial $f(x,z)$ into
irreducible factors over $\CC$, and determine which component is
$X_0(M)$.  Let $f_0(x,z)$ be the polynomial defining $X_0(M)$.  It is
easy to check that a character in $X(M)$ comes from a reducible
representation if and only if $x = z$.  So $X_0(M)$ contains a
non-integral reducible representation if and only if the polynomial
$g(x) = f_0(x,x)$ has a non-integral root.

So the hard part of checking whether Theorem~\ref{main_thm} applies is
determining the factor $f_0$ of $f$.  First, since $f(x,z)$ has
rational coefficients, there is an algorithm for factoring it over
$\CC$.  This is because one can a priori determine a finite extension
$k$ of $\QQ$ such that the irreducible factors of $f$ over $k$ are the
same as those of $f$ over $\CC$.  Take a rational line $L$ in $\CC^2$
which has simple intersections with the algebraic set $V = \{f = 0\}$
and such that $L$ and $V$ don't intersect at infinity in $\CP^2$.
Then take $k$ to be $\QQ$ adjoin the coordinates of $L \cap V$.  The
factoring of a multivariable polynomial with coefficients in a number
field is a well-studied problem (for surveys see \cite{Kaltofen82,
  Kaltofen90}).  It is worth mentioning that $f$ sometimes factors
into more components over $\CC$ than over $\QQ$, the lexicographically
smallest example being $K(45,29)$.  The computations for
Theorem~\ref{2-bridge} were done using the computer algebra system
\texttt{Maple} \cite{Maple6}, which has a built in procedure for
factoring polynomials over $\CC$.

For most of the 2-bridge knots included in Theorem~\ref{2-bridge},
every factor $f_i$ of $f$ contained a non-integral reducible
representation, and so it was not necessary to determine which $f_i$
defined $X_0(M)$.  In the exceptional cases $\{(15, 11)$, $(27, 5)$,
$(27, 11)$, $(27, 17)$, $(27, 19)$, $(33, 23)$, $(35, 29)\}$, we used
the result of Section~6.4 of \cite{HildenLozanoMontesinos96}, who
determined $f_i$ for $p < 40$ (\cite{HildenLozanoMontesinos96} gives
an algorithm for determining $f_0$ in general, but it is quite
involved).

Finally, to apply Theorem~\ref{main_thm} we have to check that $M$ is
generic.  Reid \cite[\S 4]{Reid91} showed that the only arithmetic
knot complement in $S^3$ is the figure-8 complement, which fibers. So
$M$ is non-arithmetic.  We also need to check that the cusp of the
commensurator is non-rigid.  We did this by checking that the cusp
shape is not in $\QQ(i)$ or $\QQ(\sqrt{-3})$.  Let $\Gamma \subset
\PSL{2}{\CC}$ be the image of the discrete faithful representation of
$\pi_1(M)$.  Conjugate $\Gamma$ so that the meridian generators are:

\[
a = \left(\begin{array}{cc} 1 & 1 \\ 0 & 1 \end{array}\right) \mtext{and}
b = \left(\begin{array}{cc} 1 & 0 \\ u & 1 \end{array}\right)
\] 
for some $u \in \CC$.  Riley showed that $u$ is always an algebraic
integer \cite[\S 3]{Riley72}.  Thus $\Gamma$ consists solely of
matrices with algebraic integer entries.  A longitude in the same copy
of $\pi_1(\bdry M)$ as $a$ has the form

\[\left(\begin{array}{cc} 1 & \tau \\ 
    0 & 1 \end{array}\right).
\] 
The cusp shape of $M$ is $\tau$, and so the cusp shape of $M$ is
always an algebraic \emph{integer}.  The integers in $\QQ(i)$ and
$\QQ(\sqrt{-3})$ are discrete, and so it's easy to check
numerically using \SnapPea\ \cite{SnapPea} that the cusp shape of $M$
is not in $\QQ(i)$ or $\QQ(\sqrt{-3})$, and hence that $M$ is generic.

It would have been nicer to prove that every non-fibered 2-bridge
knot satisfies the hypothesis of Theorem~\ref{main_thm}, but this
seems a difficult thing to do---in some cases there are
components of $X(M)$ which contain no reducible representations, and
it is hard to see any special property $X_0(M)$ that would prevent
this from happening there.

\begin{rmk}
  Actually, worrying about how $f$ splits up over $\CC$ as opposed to
  $\QQ$ is not really necessary.  The character varieties and all the
  maps between them in the proof of Theorem~\ref{main_thm} are all
  defined over $\QQ$.  Thus we could weaken the hypothesis of having a
  non-integral reducible representation in $X_0(M_1)$ to having a such
  a representation in the $\QQ$-irreducible component of $X(M)$
  containing $X_0(M_1)$.
\end{rmk}

\section{Examples of fibered and non-fibered pairs}\label{fiber_and_non}

\subsection{The dodecahedral knots}
The two dodecahedral knots $D_f$ and $D_s$ were introduced by
Aitchison and Rubinstein in \cite{AitchisonRubinstein92}.  They are a
pair of knots in $S^3$.  The complements $M_f$ and $M_s$ are
hyperbolic and both are quotients of $\HH^3$ by subgroups of the
symmetry group of the tiling of $\HH^3$ by $\{5, 3, 6\}$-ideal
dodecahedra.  Thus $M_f$ and $M_s$ are commensurable.  In
\cite{AitchisonRubinstein92}, they show that $M_f$ is fibered.  On the
other hand, $M_s$ is not fibered because its
Alexander polynomial is non-monic:
\[
\Delta_{M_s}=25t^4 - 250t^3 + 1035t^2 - 2300t + 2981 - 2300t^{-1} +
1035t^{-2} - 240t^{-3} + 25t^{-4}. 
\] 
The commensurator of $M_f$ and $M_s$ has a rigid
cusp, and so $D_f$ and $D_s$ are not generic (though they are
non-arithmetic). So the hypotheses of Theorem~\ref{main_thm} are
not satisfied by $M_s$.

\subsection{Small volume examples}
Here is an example of two
1-cusped manifolds which are commensurable where one is fibered and
the other not.  The two manifolds are $M_1 = m035$ and $M_2 = m039$ from
the Callahan-Hildebrand-Weeks census \cite{CallahanHildebrandWeeks}.
These manifolds have the same volume, $3.177293278...$, and same first
homology group, $\ZZ/4\ZZ \oplus \ZZ$.  Weeks' computer program \SnapPea\ 
\cite{SnapPea}, checks that $M_1$ and $M_2$ have a common 2-fold cover.
Presentations for the manifolds' fundamental groups are:
\[
\pi_1(M_1) = \spandef{a,b}{a b^3 a^{-2} b^3 a b^{-2}} \mtext{and} 
\pi_1(M_2) = \spandef{a,b}{a b^4 a b^{-1} a^{-1} b^2 a^{-1} b^{-1} }.
\] 
An easy calculation shows that the Alexander polynomials are:
\[
\Delta_{M_1} = 3 t - 2 + 3 t^{-1} \mtext{and} \Delta_{M_2} = t - 6 + t^{-1}.
\] 
Because of the lead coefficient of $\Delta_{M_1}$, the manifold $M_1$ does
not fiber over $S^1$.  On the other hand, $M_2$ is the punctured torus
bundle over $S^1$ with monodromy $+L^4 R$.  

Goodman's program \Snap\ \cite{Snap}, calculates that these manifolds
are not arithmetic and that the cusp field is a cubic extension of
$\QQ$.  Thus they are generic. It is not too hard to check that there are
non-integral reducible representations on $X_0(M_1)$.
%
Thus
this example shows that the hypothesis of Theorem~\ref{main_thm}
requiring knot complements in a $\ZZ/2\ZZ$-homology sphere is
necessary.  Some further examples among the census manifolds are:
\begin{itemize}
\item The pair $(m037, m040)$ have a common 2-fold cover and $m037$
  doesn't fiber but $m040$ does.  This pair is also commensurable with
  $m035$ and $m039$ via 4-fold covers.  
  
\item The pair $(m139,m140)$ have a common 4-fold cover, and $m139$
  doesn't fiber but $m140$ does.  Both of these manifold are
  arithmetic and so not generic.
\end{itemize}
None of these examples are knot complements in $\ZZ/2\ZZ$-homology spheres.

The strategy for finding these examples was this.  First, we used the
data provided with \Snap\ to get a list
of census manifolds grouped by commensurability invariants.  Then we
used Lackenby's taut ideal triangulations \cite{Lackenby2000} to
identify many census manifolds which fiber over $S^1$.  From this,
we selected pairs of manifolds whose trace field and cusp density were
the same, one of which fibered and the other of which did not appear
to fiber.  Most of the census manifolds fiber, making examples rare.

\subsection{Surgeries on the Whitehead link}
Let $W$ be the complement of the Whitehead link in $S^3$.  Let
$W(a,b)$ denote the 1-cusped manifold obtained by filling in one of
the two cusps of $W$ via $(a,b)$ Dehn filling.  Hodgson, Meyerhoff, and
Weeks gave a very elegant construction of a family of
fibered/non-fibered pairs which are fillings of $W$
\cite{HodgsonMeyerhoffWeeks}.  They showed that
\begin{theorem}
  Let $m \in \ZZ$ be a  multiple of $4$, $m \not\in \{0, 4\}$.  Then $M_1 =
  W(m, -1-(m/2))$ and $M_2 = W(m, -1)$ are a pair of 1-cusped hyperbolic
  3-manifolds such that:
  \begin{itemize}
    \item $M_1$ and $M_2$ have a common 2-fold cover with two cusps.
    \item $M_1$ does not fiber over $S^1$ because its Alexander polynomial 
     is not monic.
    \item $M_2$ fibers over $S^1$, being the punctured torus bundle
      with monodromy $\pm R L^m$.
\end{itemize} 
\end{theorem}

These examples overlap with the ones in the preceding section.
Namely, the pairs $(m035$, $m039)$ and $(m037, m040)$ are of this type.
The manifold $M_1$ doesn't satisfy the hypotheses of
Theorem~\ref{main_thm} because $H_1(M_1; \ZZ) = \ZZ \oplus \ZZ/m\ZZ$.  As
$m$ is divisible by $4$, $M_1$ is not a knot complement in a $\ZZ/2\ZZ$
homology sphere.

%
%

\providecommand{\bysame}{\leavevmode\hbox to3em{\hrulefill}\thinspace}
\renewcommand{\MR}{\relax\ifhmode\unskip\space\fi MR }

\end{document}